\documentclass[12pt,leqno]{article}
\usepackage{amstext}
\usepackage[dvips]{graphics}
\usepackage{theorem}
\usepackage{amssymb}
\usepackage{latexsym}
\usepackage{QED}
\usepackage{psfrag}
\input{diagrams}

\newtheorem{thm}{Theorem.}[section]
\newtheorem{prop}[thm]{Proposition.}
\newtheorem{cor}[thm]{Corollary.}
\newtheorem{lem}[thm]{Lemma.}
\newtheorem{rem}[thm]{Remark.}

\newtheorem{defn}[thm]{Definition.}
\newtheorem{exmp}[thm]{Example}
\newtheorem{prop-def}[thm]{Proposition-Definition.}

\def\pmb#1{\setbox0=\hbox{#1}%
\kern-.025em\copy0\kern-\wd0
\kern.05em\copy0\kern-\wd0
\kern-.025em\raise.0433em\box0}

\newcommand{\lra}{\longrightarrow}

\newcommand{\be}{\begin{enumerate}}
\newcommand{\ee}{\end{enumerate}}

\newcommand{\br}{\begin{array}}
\newcommand{\er}{\end{array}}

\newcommand{\gchi}{\raisebox{.4ex}{$\chi$}}

\newcommand{\up}{^{\prime}}

\newcommand{\stru}{ ^{\star}}
\newcommand{\strd}{ _{\star}}
\newcommand{\la}{\langle}
\newcommand{\ra}{\rangle}
\newfont{\ninemsbm}{msbm10 scaled 0900}
\newfont{\tenmsbm}{msbm10 scaled 1200}
\newfont{\nineeufb}{eufb10 scaled 0900}
\newfont{\teneufb}{eufb10 scaled 1200}
\newfont{\teneusm}{eusm10 scaled 1200}
\newfont{\nineeusm}{eusm10 scaled 0900}
\newfont{\bal}{cmmib10 scaled 0900}
\newfont{\tencmmib}{cmmib10 scaled 1000}

\newcommand{\sbm}[1]{\mbox{{\tenmsbm {#1}}}}

\newcommand{\fasm}[1]{\mbox{{\teneusm {#1}}}}

\newcommand{\bm}[1]{\mbox{{\boldmath ${#1}$}}}

\begin{document}
\begin{center}{\Large\bf A Generalised Harbourne-Hirschowitz
    Conjecture}\vspace{2.5ex}\end{center}
\begin{center}J. E. Alexander\end{center}

\begin{abstract}We give a generalised Harbourne-Hirschowitz
conjecture which suggests a test for determining when a
linear system on a generic rational surface separates $k$-clusters. In
particular when it is base point free or very ample.\footnote{Mathematical rewiews subject classification 14C20, 14E25, 14J26}\end{abstract}

\section{Introduction}
Recall the general context in which the Harbourne-Hirschowitz conjecture
is formulated. This is best presented as a conjecture for generic
rational surfaces which we define as follows. Fix a field $k$ and for
an integer $r\geq 0$, let $\Lambda_r$ be the residual field at the
generic point of $(\sbm P_k^2)^r$. We then have the blowing-up $\pi_r
: X_r\lra \sbm P^2_{\Lambda_r}$ in the tautological sequence of points
$x_1,\ldots ,x_r$ in the projective plane $\sbm P_{\Lambda_r}^2$. One says that
$X_r$ is the generic rational surface of rank $r$, or the blowing-up
of $\sbm P^2_k$ in $r$ generic points. 

As usual, a curve $E$ on $X_r$ with $E\simeq \sbm P^1$ and $E ^2=-1$
is called an exceptional curve on $X_r$. Since $X_r$ is rational, for
any divisor class $H$ on $X_r$, either $h^0(X_r,\fasm O(H))=0$ or
$h^2(X_r,\fasm O(H))=0$. One says that an effective divisor class $H$
is {\em non-special} if $h^1(X_r,\fasm O(H))=0$.

In slightly different forms, \cite{Har3} Harbourne and \cite{Hir1} Hirschowitz 
conjectured that (see \cite{AHi}, \cite{CM1}, \cite{CM2}, \cite{GLS} for work on the conjecture)

\noindent{\bf Conjecture (H-H)} (Harbourne-Hirschowitz) 
{\em The effective divisors $H$ on $X_r$ that are  non-special are
  those that satisify $H.E\geq -1$ for all exceptional curves $E$ on
  $X_r$.\vspace{.5ex}}

This can also be stated in the following form : if $H$ is an effective,
special divisor on $X_r$ then $H=2E+H\up$ where $H\up.E\leq 0$. In
particular $H$ is non-reduced.

\noindent{\bf Remark.} More recently, it has been pointed out by
  Ciliberto and Miranda that this conjecture is in fact equivalent to an
  older conjecture of Segre which says that a special effective
  divisor is non-reduced. 

Our purpose here is to propose a conjecture in the same vein for
separation properties such as base point freeness and very ampleness.
In \cite{D'A-H}, when studying certain classes of very ample and base point
free divisors, the authors tried unsuccessfully to propose a
suitable conjecture. In \cite{Har2}, \cite{Har3} and \cite{Di Roc}, the
respective authors showed that in the case $r\leq 9$, base point
freeness, very ampleness and more general separation properties were
characterised by intersection numbers on a small class of curves, but
no general conjecture was forthcoming. 

In accordance with known results for $r\leq 9$ the following
conjecture says that if a divisor class $H$ with sufficiently many
sections does not separate $k$-clusters then there is an integral
(even smooth and irreducible) curve $E$ such that $H$ does not
separate some $k$-cluster in $E$. Furthermore, evidence would suggest
that this failure must be detectable numerically. Before formulating
the conjecture we need the following  
\begin{defn} {\em A $k$-{\em cluster} $Z$ on $X_r$ is any finite closed subscheme of $X$ of
length (i.e. degree) $k$. We say that an effective divisor class $H$
on $X$ separates $k$-clusters if $h^0(X_r, I_Z(H))=0$.}\end{defn}

\noindent{\bf Conjecture 1.} {\em Let $k\geq 0$ be an integer. If $H$ is an effective divisor class on $X_r$ such that $\gchi (X_r,\fasm O(H))\geq 3k$, then $H$ separates $k$-clusters if the following (necessary) conditions are satisfied : for all integral
  curves $E$ of genus $a$, if $0\leq a\leq k$ then $H.E\geq 2a-1+k$
  and if $k< a \leq \frac{4k}{3}$ then $H.E\geq a+2k-1$.}

A weakness in this formulation of the conjecture is that it would
appear to be necessary to test on all integral classes of low
genus. In fact for $r\geq 3$ one can use tha Weyl group action on
$X_r$ to obtain a {\em standard} representative of classes $H$ that
satisfy $H.E\geq 0$ for all exceptional classes $E$. One can then
show that the conjecture is equivalent to testing on a finite set of
isolated, integral, standard classes of low genus.

In \S2, we give the structure theorem for effective divisors on $X_r$
as implied by H-H. In \S3, we give the motivation and assumptions
underlying the conjecture. 
In \S4, we review the Weyl group action and the notion of a standard
class. While the notion is not new, I have not seen mention in the
literature of the caracterisation of these classes as minimising for
the intersection product in all semi-standard orbits
(prop. \ref{min}).This observation allows the reformulation of the
conjecture given in \S4.3. 

In the remainder of the introduction we make some preliminary
remarks.
\subsection{Preliminary remarks}
\be 
\item If $k=0$, conjecture 1 is just the Harbourne-Hirshowitz conjecture for
  the natural cohomologie of $H$. For $k>0$, it corresponds to the
  surjectivity of the canonical map
\begin{equation}\label{rest}H^0(X_r,\fasm O(H))\lra H^0(X_r,\fasm O_Z(H)) \end{equation}
for all $k$-clusters, so that it is equivalent to $|H|$ being base point
free (resp. very ample) when $k=1$ (resp. $k=2$). Note that in the
latter cases ($k=1,2$) we only require that $H.E\leq 2a-1+k$ for
$0\leq a\leq k$.
\item The conjecture is true for $r\leq 9$. This is best seen using
  the equivalence of conjectures 1 and 2.  When $k=0,1,2$ this was
  proven in \cite{Har2} Harbourne  and, for $k\geq 3$, $r\leq 8$ in
  \cite{Di Roc} Di Rocco, but any reader will see for themselves that
  her proof using the general vanishing theorem [S] goes over to $r=9$.
\item The condition $\gchi (X_r,\fasm O(H))$ is the natural condition
  for a {\em general}  $H$ to separate $k$-clusters in the following
  sense. Consider the canonical diagramme
$$\begin{diagram}[size=2em,postscript]
Z&\rInto&\mbox{Hilb}^k(X_r/\Lambda)\times_{\Lambda}X_r&\rTo ^{p}&X_r\\
&&\dTo_q&&\\
&& \mbox{Hilb}^k(X_r/\Lambda)&&\end{diagram}$$
where $Z$ is the tautological subscheme. We then have the canonical
vector bundle map
$$H^0(X_r,\fasm O(H))\otimes\fasm O_{\mbox{Hilb}^k(X_r)}\lra q\strd \fasm O_Z(p\stru H) $$
giving the map (\ref{rest}) at each point of the Hilbert scheme. The
left hand side has rank $n=h^0(X_r,\fasm O(H))$ and the right hand
side has rank $k$, so that the degeneracy locus is empty or has
codimension $\leq n-k+1$. If the map is to be nowhere degenerate in
general (i.e. without further numerical conditions being applied) then
we need $n-k+1>\mbox{dim Hilb}^k(X_r)=2k$, giving the condition
$h^0(X_r,\fasm O(H))\geq 3k$. In view of the Harbourne-Hirschowitz
conjecture, we should have $h^1(X_r,\fasm O(H))=0$ and since
$h^2(X_r,\fasm O(H))=0$ for an effective divisor on the rational
surface $X_r$, this should be the same thing as $\gchi (X_r,\fasm
O(H))\geq 3k$.
\ee
\section{About isolated curves}
\subsection{Harbourne-Hirschowitz and the structure of effective classes}
\begin{defn}\label{isolated}{\em  An effective divisor $E$ on $X_r$
    satisfying $E ^2 = a-1=E.K$; $a\geq 0$; will be called an isolated
    curve of genus $a$. {\em Note that exceptional curves are just
    isolated curves of genus zero}.}\end{defn} 

If one supposes that H-H holds, one gets a precise description of the
effective classes. In particular, one gets that for $a\geq 2$, an
isolated curve of genus $a$ is reduced and irreducible, and isolated in
its linear system.
\begin{prop}(Bossini)\label{str} If H-H holds, then for any effective
 divisor class $H$ on  $X_r$, the generic curve $C$ in $|H|$ can be
 decomposed as an orthogonal sum $A+n_1F_1+\cdots +n_rF_r$, where each
 $F_i$ is an exceptional curve, $n_i>0$ and  $A$ is effective and
 satisfies $p(A)\geq 0$, $A^2\geq p(C)-1$ and $A ^2\geq 0$ if
 $p(A)=0$; $p$ being the arithmetic genus. If $A ^2>0$ then $A$ is
 reduced and irreducible. If $A ^2=0$ then we have one of the
 following two cases 
\be 
\item $A ^2=0$, $A.K=-2b$ ($b>0$) in which case $A=bE$, where $E$ is a
  pencil class (definition \ref{classes})
\item $A^2=0=A.K$ in which case $A=mE$ ($m>0$) where $E$ is a smooth
  irreducible elliptic classe which, with respect to a suitable exceptional
  configuration (see \S~\ref{excep}), can be written $E\equiv 3E_0-E_1-\cdots -E_9$.
\ee
\end{prop}

\begin{proof}$\!\!\!\!\!.$  One can write $C=A+\sum_{i=1}^s n_iE_i$, where
  the $E_i$ are the reduced, irreducible components of $C$ which are
  rational of self-intersection $<0$. . Now, it is known that the only 
  integral rational curves on $X_r$ having self-intersection $<0$
  are the exceptional curves. In fact if $W$ is such a curve with
  $W^2\leq 2$, then $W.(-K)<0$ and,
  specialising to the situation where all points are generic on
  a smooth plane cubic, $W$ must specialise to an effective divisor
  having a smooth elliptic component, but this is impossible for $W$
  rational. We conclude that each $E_i$ is exceptional.

Since each $E_i$ is isolated we have $E_i.E_j=0$ otherwise
$\mbox{dim}|E_i+E_j|>0$ and the generic curve in $|E_i+E_j|$ does not
decompose. As well $A$ is non-special by H-H, because it is positive
on all exceptional curves.
The isolation of $\sum_{i=1}^s n_iE_i$ in $|H|$ gives $h^0(X_r,\fasm
O(A))=h^0(X_r,\fasm O(C))=h^0(X_r,\fasm O(A+E_i))$ which shows
that $h^0(E_i,\fasm O_{E_i}(A+E_i))=0$ and finally that $A.E_i\leq 0$,
but $A.E_i\geq 0$ in any case.

Now if $A=A_1+A_2$ where $A_2$ is a reduced and irreducible curve and
$A_1$ is an effective divisor, then $A_2$, not being 
rational of negative self-intersection, is  non-special by
H-H. This implies that
$A_2^2\geq p(A_2)-1\geq 0$ if $p(A_2)>0$ and $A_2^2\geq 0$ if
$p(A_2)=0$, where $p$ denotes the arithmetic genus. If $A_2^2=0$ then
$A_2$ is a pencil class or an isolated integral elliptic class. 

Since this is true for all integral components $A_2$ of $A$, we
conclude that $A ^2\geq 0$ and even $A ^2\geq p(A)-1$ (using the
standard formulas for $p$). As such $A ^2\geq 0$ and $A ^2>0$ unless
$A$ is a pencil class or an isolated elliptic class. 

Now using the fact that $A_1$ and $A_2$ are 
non-special (H-H), and the surjectivity of the map
$|A_1|\times|A_2|\lra|A_1+A_2|$, we find by H-H 
$$\br{lll}
\mbox{dim}\, |A_1|+\mbox{dim}\, |A_2|&\geq &\mbox{dim}\, |A_1+A_2|\\
                                     &=&\gchi (X_r,\fasm
                                     O(A_1+A_2))-1\\
                                     &=&\gchi (X_r,\fasm
                                     O(A_1))-1 +\gchi (X_r,\fasm
                                     O(A_2))-1 +A_1A_2\\
                                     &=&\mbox{dim}\,
                                     |A_1|+\mbox{dim}\,
                                     |A_2|+A_1.A_2\\
\er$$
so that $A_1.A_2=0$.

By the algebraic index theorem, since we know that $A_i ^2\geq 0$,
$i=1,2$, we must have $A_1^2=A_2^2=A_1.A_2=0$ and $C=bA_2$ (by the
Cauchy-Schwartz inequality after writing $A_1$ and $A_2$ in an
exceptional configuration) with $A_2$ a pencil class, or
$A_2$ an  isolated integral elliptic class, hence equivalent to
$3E_0-E_1-\cdots -E_9$ by lemma \ref{bossi}.\end{proof}

It is conjectured that all isolated curves $E$ of genus $a\geq 2$ are
in fact smooth as well (see remark~\ref{ice}). 
\subsection{The ample classes}
\begin{prop} If H-H holds then a divisor class $H$ on $X_r$ is ample
  if and only if $H^2>0$ and $H.E>0$ for all exceptional curves $E$ on
  $X_r$. If $H$ is standard and $H ^2>0$, and  $H.E\leq 0$ for some
  reduced and irreducible curve $E$ on $X_r$, then $E$ is exceptional
  and $H.E=0$.\end{prop}

\begin{proof}$\!\!\!\!\!.$ If $H$ is
  standard then $h^2(X_r,\fasm O(H))=0$ and the condition $H^2>0$
  shows that $aH$ is effective for $a>\!\!>0$. We can thus suppose $H$
  effective, hence reduced and irreducible by (\ref{str}). If $E$ is a
  reduced and irreducible curve then $H.E\geq 0$, but if $H.E=0$ then
  $E^2 < 0$ by the algebraic index theorem and by (\ref{str}) again, this
  implies that $E$ is exceptional. \end{proof}
In this direction one might consult \cite{Xu}.

\section{The conditions of the conjecture}

\subsection{The conditions $H.E\geq 2a-1+k$, $a\leq k$ and $H.E\geq
  a-1+2k$, $k<a\leq \frac{4}{3}k$ of the conjecture} 
Throughout this section let $H$ be an effective divisor class on $X_r$
($r\geq 3$) such that 
$\gchi (\fasm O(H))\geq 3k$ and $H.F\geq k-1$ for all exceptional
curves $F$ sur $X_r$.

It has already been said that the underlying assumption in conjecture
1. is that if $|H|$ does not separate $k$-clusters then this because
it does not separate some $k$-cluster in an integral (even smooth)
curve on $X_r$. It is also supposed that if the map 
\begin{equation}\label{reste}\rho : H^0(X_r,\fasm O(H))\lra H^0(X_r,\fasm
  O_E(H))\end{equation}
has rank $\geq 2k$ then $|H|$ separates $k$-clusters in $E$. We will
see that if H-H holds and $\fasm O_E(H)$ is special, or if $H-E$ is
effective and special, then (\ref{reste})
always has rank $\geq 2k$. In the remaining cases (\ref{reste}) is
surjective and there are two possibilities. Firstly, if $k>a$ then it
is well known that
$\fasm O_E(H)$ separates $k$-clusters if and only if $H.E\geq
2a-1+2k$, but in this case $h^0(\fasm O_E(H))$ can be less than
$2k$. Secondly, if $k\leq a$, a necessary condition for $\fasm O_E(H)$
to separate $k$-clusters is that $H.E\geq a-2+2k$, and this is
sufficient if $\fasm O_E(H)$ is a general invertible sheaf of degree
$H.E$. We are thus supposing in the conjecture that $\fasm O_E(H)$ is
general in this sense.

This said we must prove the
\begin{prop} Suppose that H-H holds. Then the image of
  (\ref{reste}) has dimension $<2k$ only if $\fasm O_E(H)$ is
 non-special. When $\fasm O_E(H)$ is non-special, $a>k$ and $H.E\leq a-2+2k$
 one has $4k\geq 3a$.\end{prop}

\begin{proof}$\!\!\!\!\!.$ Firstly, if H-H holds, $E$ is non-special.

There are two more or less obvious cases where the image of
(\ref{reste}) has dimension $\geq 2k$. These are when $H-E$ is not
effective (obvious) and when $H-E$ is effective and special. In fact,
in  the second case, if H-H holds, then by (\ref{str}) we
can write $H-E=A+\sum_i n_iF_i$ as an orthogonal sum of an effective
standard divisor $A$ and multiples $n_i>0$ of exceptional curves $F_i$
with at least one of the $n_i\geq 2$. Since $k-1\leq H.F_i=E.F_i-n_i$ one finds
$E.F_i\geq n_i+k-1$ and 
$$\br{lll}
&    &    h^0(X_r,\fasm O(H))-h^0(X_r,\fasm O(H-E))\vspace{.5ex}\\
&=   &    \overline{\gchi} (H)-\overline{\gchi} (A)\vspace{.5ex}\\
&=   &    \overline{\gchi} (A)+\overline{\gchi} (E+\sum_i n_iF_i)+A.E-
          \overline{\gchi} (A)\vspace{.5ex}\\
&=   &    -\sum_i n_i(n_i-1)/2+\sum_i n_iE. F_i +A.E\vspace{.5ex}\\
&\geq&    \sum_i n_i(n_i-1)/2+k( \sum_i n_i)+A.E\vspace{.5ex}\\
&\geq   &    2k
\er$$
as required. 

Moreover, since $H.F\geq 0$ for all exceptional curves $F$ on $X_r$,
$H$ is non-special. Also if $\fasm O_E(H)$ 
is special, $h^2(\fasm O_E(H-E))>0$, implying that $H-E$ is not
effective, so that the image of (\ref{reste}) is at least $2k$
dimensionl in this case as well.

We can thus suppose henceforth that $H-E$ and $\fasm O_E(H)$ are
effective and non-special and we will show that if $a>k$ and $H.E\leq
a-2+2k$ one has $4k\geq 3a$.

\noindent{\bf Claim.} In this case, $H-3E$ is not effective, $\fasm
O_E(H-E)$ is non-special, and $\fasm O_E(H-2E)$ is special. 

This gives the
desired inequality as follows using Riemann-Roch and Cliffird's
theorem
$$\br{lll}\frac{1}{2}
(E.(E-2H))+1&\geq & h^0(\fasm O_E(H-2E))\\
                &\geq & h^0(\fasm O(H-2E))  \\
                &\geq & h^0(\fasm O(H))- h^0(\fasm O_E(H))- h^0(\fasm
                                                           O_E(H-E))\\
                &=&     h^0(\fasm O(H))-2E.H+E ^2 +2(a-1)\er$$
so that
$$5E.H +2\geq 2h^0(\fasm O(H))+4E ^2+4(a-1)$$
and using $E ^2\geq a-1$ (\ref{st}) and $a-2+2k\geq H.E$, we get $4k\geq
3a$.

\noindent{\bf Proof of the claim.} Since $(H-3E).E = H.E-3E ^2 \leq a-2+2k-3(a-1)=2(k-a)+1<0$ and $E
^2\geq 0$ it follows that $H-3E$ is not effective. If $\fasm O_E(H-E)$ is
special, then $h^2(\fasm O(H-2E))=h^1(\fasm O_E(H-E))\neq 0$, so that
$H-2E$ is not effective and by Clifford's theorem 
$$k+1\leq h^0(\fasm O(H-E))\leq h^0(\fasm
O_E(H-E))\leq\frac{1}{2}(E.H-E ^2)+1\leq\frac{1}{2} (a-2+2k-a+1)+1=
k+\frac{1}{2}$$
From this contradiction we conclud that $\fasm O_E(H-E)$ is
non-special. Finally,  
$$\br{lll}\overline{\gchi} (\fasm O(H-3E))&=& 
               \overline{\gchi} (\fasm O(H))-\overline{\gchi}(\fasm
               O(3E))-3E.(H-3E)\\
 &=& \overline{\gchi} (\fasm O(H))-3\overline{\gchi}(\fasm O(E))+6E
               ^2-3E.H\\
 &\geq & 3k+6(a-1)-3(a-2+2k)\\
&=& 3(a-k) \; >\; 0
\er$$
so that $h^2(\fasm O(H-3E))>0$ and $\fasm O_E(H-2E)$ is special.
\end{proof}

\section{\label{excep}Translation into standard divisors}
As we will see below, standard divisors are essentially those that
test positive on all exceptional curves and include all
integral classes that are not exceptional. The more general
semi-standard classes include all efective classes on an $X_r$. What
we want here is threefold: (1) an intrinsic  
characterisation of the (semi)-standard classes, (2) a standard way of
writing such a class with respect to a suitable exceptional
configuration, (3) the unicity of this expression (but not of the
exceptional configuration). I have found no reference in the
literature to parts (1) and (3), but (2) has been used extensively
since Noether. We thus begin by defining the classic $\bm
E$-(semi)-standard classes for a fixed exceptional configuration $\bm
E$ (definition~\ref{standard}), and then show that this is equivalent
to an intrinsic definition (proposition~\ref{st}). Finally in
proposition~\ref{es}, we deal with the unicity via
a minimising characterisation of the exceptional configurations that give
$\bm E$-(semi)-standard representations of  (semi)-standard
classes. This gives corollary~\ref{min} which is essential for the
reformulation of the conjecture.

Before giving the definitions
let us recall \cite{Hir2} that the Weyl groupe $W_r$ of the surface $X_r$
(see \cite{Har3} and \cite{Dol}) can be viewed as a groupe of
$k$-automorphismes of $X_r$. There is an induced action of $W_r$ on $\mbox{Pic}(X_r)$ which
stabilises the intersection form and fixes the canonical class. The
latter action is faithful so that $W_r$ is identified with a subgroup
of the orthogonal group $O_r$ of $(\mbox{Pic}(X_r), \la\; ,\; \ra)$
where $\la\; ,\; \ra$ is the intersection form.

This induced action on $\mbox{Pic}(X_r)$ is simply transitive on the
{\em exceptional configurations} which are by definition the orthogonal
sequences of effective, irreducible classes
$$\bm E=(E_0,E_1,\ldots ,E_r)$$
satisfying
$$E_0 ^2=1\quad , \quad E_0.K_r=-3\quad , \quad E_j ^2=-1=E.K_r\quad  ,\quad  E_i.E_j=0\quad  , \quad 0\leq i<j\leq r$$
In particulier $E_i\simeq \sbm P^1$ for $i=0,1,\ldots ,r$. For any
fixed exceptional configuration $\bm E$ as above, the group $W_r$  is
generated by the orthogonal reflexions $\sigma_i=\mbox{id}+\la r_i,\;
\ra r_i$ where 
$$r_0=E_0-E_1-E_2-E_3,\hspace{1ex} r_i=E_i-E_{i+1},\hspace{1ex} i=1,\ldots ,r-1$$
Central to the proof of these results is Noether's inequality (see \cite{Dol}
§ 5), which
says that if $r\geq 3$, then a  sequence of non-negative integers
$d,m_1\geq m_2\geq \cdots \geq m_r$ satisfying 
$$d^2-m_1^2-\cdots -m_r^2=-1=-3d+m_1+\cdots +m_r$$ 
also satisfies $d< m_1+m_2+m_3$. This inequality
can also be used to show that for $r\geq 3$ any irreducible rational
class $E$ (i.e. $E$ is irreducible, effective and $E ^2=-2-E.K_r\geq
-1$ lies in the $W_r$ orbit of one and one only of the following classes
\begin{equation}\label{rc}E_r,\hspace{.25ex}E_0-E_1,\hspace{.25ex}E_0,\hspace{.25ex}2E_0,\hspace{.25ex} dE_0-(d-1)E_1-E_2,\hspace{.25ex} dE_0-(d-1)E_1,\quad d\geq 2\end{equation}
Note that these are uniquely determined by $E ^2$ except for the pair
$2E_0$ and $3E_0-2E_1-E_2$.
\begin{defn}\label{classes}{\em  An irreducible rational classe $E$ on $X_r$ with $E ^2=0$
(resp. $E ^2=1$, resp. $E ^2=2$) will be called a pencil class (resp. line
class, resp. quadratic class). Note that this is $W_r$-equivalent to
$E_0-E_1$ (resp. $E_0$, resp. $2E_0-E_1-E_2$).}
\end{defn}
It is also true that $W_r$ acts transitively on the orthogonal
sequences of exceptional curves of length $\neq r-1$. When the length
is $r-1$ there are two orbits which, for a fixed exceptional
configuration $\bm E=(E_0,E_1,\ldots ,E_r)$, have representatives of
the form $(E_2,\ldots ,E_r)$ et $(E_0-E_1-E_2,E_3,\ldots ,E_r)$. These
correspond respectively to sequences of length $r-1$ that can be
extended and those that cannot. If we forget the ordering, then
the non-extendable orthogonal sequences of exceptional curves
correspondent to quadratic classes $E$ which induce birational
morphismes of $X_r$ to a smooth quadrique $Q$ in $\sbm P ^3$,
contracting the exceptional curves in such a non-extendable sequence.

We define standard classes as follows
\begin{defn}\label{standard} {\em A divisor class $H$ on $X_r$ is said to be $\bm
 E$-{\em standard} for an exceptional configuration $\bm E=(E_0,E_1,\ldots
 ,E_r)$ if 
$$H\equiv dE_0-m_1E_1-\cdots -m_rE_r$$
where $d\geq m_1\geq \cdots \geq m_r\geq 0$, $d\geq m_1+m_2$ and
$d\geq m_1+m_2+m_3$. 

We say that $H$ is $\bm E$-{\em semi-standard} if
$d\geq 0$, $d\geq m_1+m_2+m_3$ et $d\geq m_1\geq m_2\geq \cdots \geq
m_r$ (but some $m_i$ may be $<0$).

We will say that a divisor class $H$ on $X_r$ is {\em standard}
(resp. {\em semi-standard}) if it is $\bm
E$-standard (resp. $\bm E$-semi-standard) for some exceptional configuration $\bm E$ on $X_r$.}
\end{defn}
\begin{exmp} The class $5E_0-2E_1-2E_2-2E_3$ is not $\bm
  E$-standard, but is standard as one can see by making the standard
  quadratic transformation by the reflection with root $r_0$.\end{exmp}
\subsection{generating classes}
It is obvious and was pointed out by Harbourne that the $\bm E$-standard
classes are precisely those that can be expressed as non-negative sums
\begin{equation}\label{gcp}H\equiv aE_0 + b(E_0-E_2)+b(2E_0-E_1-E_2)+\sum_{i=3}^r\alpha_iC_i\end{equation} 
where $a,b,c,\alpha_i\geq 0$ and $C_i=-K_r+E_{i+1}+\cdots +E_r$
($i\geq 3$).
\begin{defn}{\em  We say that the classes $E_0$, $E_0-E_1$, $2E_0-E_1-E_2$,
$C_i$ ($i\geq 3$) are the (standard) {\em generating classes} of the
exceptional configuration $\bm E$. The first three
generating classes are called the rational generating classes while
the $C_i$, which have arithmetic genus one, are called the elliptic
generating classes.}\end{defn}
\begin{rem} {\em As we saw above the rational generating classes of self
  intersection 0,1,2 each form a single orbit under the Weyl group
  action. The same is true for the elliptic generating classes for
  each self intersection number $\leq 6$ since $W_r$ fixes $K_r$ and
  acts transitively on the orthogonal sequences of exceptional curves
  of length $\neq r-1$.} \end{rem}
\begin{lem}\label{bossi}(Bossini) If $C$ is an  $\bm
  E$-standard class satisfying $C^2=0=C.K_r$, then $C\equiv \alpha_9C_9$.  \end{lem}

\begin{proof}$\!\!\!\!\!.$ The class $C$ is effective and it is clear that $C.E>0$ on all rational
  generating classes and $C.C_i\geq 0$ since $C.C_r=C.-K_r=0$. As well
  $C.C_i>0$ for $i\leq 8$ since any effective divisor  that is a
  proper component of a curve in $|C_i|$ is rational. As such,
  $C\equiv \sum_{i\geq 9} \alpha_i C_i$. If $C\neq \alpha_9C_9$ then
  $C^2<0$.\end{proof} 
\subsection{semi-standard classes}
An $\bm E$-semi-standard class $H$ has one of two forms reflecting the
two orbits of orthogonal sequences of exceptional curves of length
$r-1$. In fact, either $H=(dE_0-m_1E_1-\cdots -m_sE_s)+n_{s+1}E_{s+1}+\cdots +n_rE_r$ where the first part is
standard and $n_i=-m_i>0$, or $m_3<0$, $m_1,m_2\geq 0$, $d <  m_1+m_2$
and 
$$\br{lll}
H&\equiv&(2d-m_1-m_2)E_0-(d-m_2)E_1-(d-m_1)E_2+(m_1+m_2-d)(E_0-E_1-E_2)\\
&&\; + n_3E_3+\cdots +n_rE_r\er$$
where $n_i=-m_i>0$. As such any $\bm E$-semi-standard class, can be
written in the form $A+\sum_i n_iF_i$, where $A$ is $\bm E$-standard,
$n_i>0$ and the $F_i$ are a familly of orthogonal exceptional curves.
\begin{prop} \label{st} For $r\geq 2$, a divisor class $H$ on $X_r$ is standard
  (resp. semi-standard) if and only if $H.E\geq 0$ for all exceptional
  curves $E$ on $X_r$ (resp. $H.E\geq 0$ for all line and pencil
  classes $E$ on $X_r$). 

A classe $H$ is semi-standard if and only if
  it is an orthogonal sum of the form $H=A+\sum n_iF_i$, where $A$ is
  standard, $n_i>0$, and each $F_i$ is an exceptional curve.
\end{prop}
\begin{proof}$\!\!\!\!\!.$ Firstly, if $H$ is non-negative on all exceptional
  classes, then choosing any line class $E_0$, we obtain an
  exceptional configuration $\bm E$, with $H=dE_0-m_1E_1-\cdots
  -m_rE_r$, $m_i\geq 0$ and $(d-m_1-m_2)=H.(E_0-E_1-E_2)\geq 0$, so
  that $H.E_0=d\geq 0$ on all line classes $E_0$. 

Suppose then that $H$ is non
  negative on all line classes. In this case we can choose a line
  class $E_0$ such that $H.E_0=d$ is the minimum for all line classes
  on $X_r$ and by suitably ordering the exceptional curves
  contracted by $E_0$ we obtain an ecxceptional configuration $\bm
  E=(E_0,E_1,\ldots ,E_r)$ with 
$$H=dE_0-m_1E_1-\cdots -m_rE_r$$
$m_1\geq \cdots \geq m_r$ and, if $r\geq 3$, 
$2d-m_1-m_2-m_3=H.(2E_0-E_1-E_2-E_3)\geq H.E_0=d$ by the minimality of
$H.E_0$ over all line classes. 

When $H$ is positive on all exceptional
curves we have $m_r\geq 0$, $d\geq m_1+m_2\geq m_1$ and when  $H$ is positive
on all pencil classes we have  $d\geq m_1$.

Conversely, if $H=dE_0-m_1E_1-\cdots -m_rE_r$ is $\bm E$-standard, we
can write 
$$H\equiv aE_0 + b(E_0-E_2)+b(2E_0-E_1-E_2)+\sum_{i=3}^r\alpha_iC_i$$
as in (\ref{gcp}) and it suffices to note that each of the
generating classes is positive on all exceptional curves.

If $H$ is $\bm E$-semi-standard, then as we saw above, $H$ is an
orthogonal sum $H\equiv A+\sum_in_iF_i$ where $A$ is $\bm E$-standard,
$n_i>0$ and the $F_i$ are exceptional. Conversely every such class is
non-negative on all line and pencil classes. 
\end{proof}
\begin{cor}\label{effst} On $X_r$ every effective class is
  semi-standard and every integral class that is not exceptional is
  standard.\end{cor} 
\begin{rem}\label{ice}{\em We have defined isolated curves of genus $a\geq 1$ to be reduced and
irreducible curves $E$ with $E ^2=a-1=E.K_r$. The preceeding
proposition says that if H-H holds then for $a\geq 2$, we could define
these to be standard classes with $E ^2=a-1=E.K_r$.

The isolated curves of genus $a\leq 4$ were completely classified in
\cite{Bos} and \cite{Mig}. The only standard
classe with $E ^2=0=E.K_r$ is $C_9$, while the only standard
classes with $E ^2=1=E.K_r$ are 
\begin{equation}\label{isol}\br{lll}G_1&=&4E_0-2E_1-E_2-\cdots -E_{12}\\
G_2&=&6E_0-2E_1-\cdots -2E_8-E_9-E_{10}-E_{11}\\
G_3&=&9E_0-3E_1-\cdots -3E_8-2E_9-2E_{10}\er\end{equation}
and it has been proven that for $a=2,3,4$ all standard classes with $E
^2=a-1=E.K_r$ contain a unique curve and it is smooth and
irreducible. }\end{rem}
\begin{prop}\label{es} Let  $H\equiv dE_0-m_1E_1-\cdots -m_rE_r$ be an $\bm
  E$-standard class, then
\be
\item for $i\geq 1$, $m_i+\cdots +m_r$ is the minimum value of
  $H.(G_i+\cdots +G_r)$ for all extendable orthogonal sequences of exceptional
  curves $G_i,\ldots ,G_r$. If $i\geq 3$ this is the minimum value on
  all sequences of exceptional curves, orthogonal or not.
\item $d-m_1-m_2+m_3+\cdots +m_{r}$ is the minimum value of
  $H.(G_2+\cdots +G_r)$ on all unextendable orthogonal sequences  of
  exceptional curves $G_2,\ldots ,G_r$.
\ee

In particular, if $H$ is a standard class which can be written in the
  $\bm E$-standard form $dE_0-m_1E_1-\cdots -m_rE_r$, then the
  sequence $d,m_1,\ldots ,m_r$ is unique, independent of the
  exceptional configuration $\bm E$ in which $H$ is $\bm
  E$-standard.
\end{prop}
\begin{proof}$\!\!\!\!\!.$ As we saw in (\ref{gcp}), we can write $H$ uniquely in the form
$$H\equiv aE_0 + b(E_0-E_2)+c(2E_0-E_1-E_2)+\sum_{i=3}^r\alpha_iC_i$$
For $i\geq 3$ we have $m_i=\alpha_i+\cdots +\alpha_r$. If $E$ is an
exceptional curve 
$$H.E\geq \alpha_i E.C_i+\cdots + \alpha_r E.C_r$$
and $E.C_i=1+E(E_{i+1}+\cdots +E_r)$, so that if $H.E < m_i$ then $E=E_j$ for some $j>i$. This shows that
for {\em all} sequences $G_i,\ldots ,G_r$ ($3\leq i\leq r$) of
exceptional curves, $H.(G_i+\cdots +G_r)\geq H.(E_i+\cdots +E_r)$. 

If $E$ is an exceptional curve and $m_3<H.E< m_2=m_3+c$ then $c>0$,
$E. (2E_0-E_1-E_2)=0$ et $E\neq E_i$ for $i=3,\ldots ,r$. As such
$E=E_0-E_1-E_2$. It follows that for extendable (resp. unextendable)
sequences $G_2, \ldots , G_r$ of orthogonal exceptional curves, the minimum value of $H.(G_2+\cdots
+G_r)$ is attained on the sequence $G_i=E_i$; $i=2,\ldots ,r$
(resp. $G_2=E_0-E_1-E_2$, $G_i=E_i$; $i=3,\ldots ,r$).

If $E$ is an exceptional curve and $m_2<H.E< m_1=m_2+b$, then $b>0$,
$E\neq E_j$; $j=2,\ldots ,r$; et $(E_0-E_1)E=0$. In this case
$E=E_0-E_1-E_j$ for some $j>1$ and $E,E_2,\ldots ,E_r$ is not an
orthogonal sequence. It follows that for orthogonal sequences of
length $r$, $G_1,\ldots ,G_r$ of exceptional curves, the minimum value of
$H.(G_1+\cdots +G_r)$ is attained on the sequence $(E_1,\ldots
,E_r)$. 
\end{proof}
\begin{cor} The decomposition of a semi-standard class $D$ as an
  orthogonal sum $A+\sum_i n_iF_i$, where $A$ is standard, the
  $F_i$ are exceptional and $n_i>0$, is unique.\end{cor}

\begin{proof}$\!\!\!\!\!.$ If $B+\sum_i m_iG_i$ is another such decomposition,
  then, $B$ being standard, $D.F_i<0$ implies that $F_i.G_j<0$ for
  some $j$. Since the $F_i$ and $G_i$ are smooth and irreducible
  classes, we can suppose that $F_i=G_i$, then $n_i=D.F_i=m_i$.\end{proof} 
\begin{prop}\label{min} Fix an exceptional configuration $\bm E=(E_0,E_1,\cdots
  ,E_r)$ on $X_r$. Let $H\equiv dE_0-m_1E_1-\cdots -m_rE_r$ be an $\bm
  E$-standard divisor and let $D$ be a semi-standard divisor on
  $X_r$. Then for $\sigma\in W_r$, the minimum value of $H.\sigma(D)$
  occurs exactly when $\sigma(D)$ is $\bm E$-semi-standard. Otherwise
  said, if  $\bm F=(F_0,F_1,\cdots ,F_r)$ is an exceptional
  configuration for which $D = d\up F_0 - m\up_1 F_1-\cdots - m\up_r F_r$ is
  $\bm F$-semi-standard and $D\up =  d\up E_0 - m\up_1 E_1-\cdots - m\up_r E_r$ then 
$$H.D\geq H.D\up$$
\end{prop}

\begin{proof}$\!\!\!\!\!.$ Write $H$ as
$$H\equiv aE_0 + b(E_0-E_2)+c(2E_0-E_1-E_2)+\sum_{i=3}^r\alpha_iC_i$$
and $D$ as an orthogonal sum (\ref{st}) $D=A+\sum_i ^r n_jG_j$ ($i\geq
2$) where $A$ is $\bm F$-standard, $n_i>0$ and $G_i,\ldots ,G_r$ is an
orthogonal sequence of exceptional curves with either $F_j=E_j$ for
$j=i,\ldots ,r$, or $i=2$, $G_2=F_0-F_1-F_2$, $G_j=F_j$ for $j=3,\ldots ,r$.

 By (\ref{es}), the proposition holds when $A=0$, so it suffices to
 prove the proposition for $D$ standard and even for $D$ one of the
 $\bm F$-standard generating classes. By (\ref{es}), the proposition holds for
 $D=-K_r-E_i-\cdots -E_r$ ($i\geq 3$) and by the symmetry of the
 proposition it suffices to suppose that $H$ (resp. $D$) is one of the
 rational generating classes.

The system $|E_0|$ (resp. $|2E_0-E_1-E_2|$) is positive on all
effective classes other than exceptional curves, so we need only show
that $(2E_0-E_1-E_2).F_0\geq 2$ and
$(2E_0-E_1-E_2).(2F_0-F_1-F_2)\geq2$. This follows from the fact that
the only effective divisors with $(2E_0-E_1-E_2).G\leq 1$ move in a
linear system of projective dimension at most one.\end{proof}
When $r\leq 1$, there is only one exceptional configuration so that
$\bm E$-standard and standard are equivalent.
\begin{rem}\label{minrat} {\em It follows from proposition~\ref{es} that if $r\geq 3$ then
  $H.E_r$ is the minimum value of $H.E$ on all exceptional curves $E$.
  It then follows from proposition~\ref{min} and the list of effective
  $\bm E$-standard rational classes~(\ref{rc}), that $H.E_r$ is the
  minimum on all effective rational classes.}\end{rem} 
\subsection{reformulation of the conjecture}
\begin{defn}\label{nd}{\em 
For an exceptional configuration $\bm E=(E_0,E_1,\ldots ,E_r)$
on $X_r$ we extend the definition of $\bm E$-standard class to
include, for any $\delta\geq 0$, those classes $E$ on $X_{r+\delta}$
which for the canonical extension
$\bm E\up=(E_0,E_1,\ldots ,E_r,E_{r+1},\ldots ,E_{r+\delta})$ of the
exceptional configuration to $X_{r+\delta}$ is $\bm E\up$ standard on
$X_{r+\delta}$. An $\bm E$-standard isolated curve of genus $a$ will
then be an integral $\bm E$-standard class (for this extended
definition)  $E$ such that $E ^2=a-1=E.K_{r+\delta}$.}\end{defn} 
\begin{rem}{\em Note that if $H$ is a semi-standard or standard class
    on $X_r$ then so is its pull back (also denoted by $H$) to
    $X_{r+\delta}$. }\end{rem} 
We can now reformulate conjecture 1. in the following form

\noindent{\bf Conjecture 2.} {\em Let $k>0$ be an integer and let $\bm
  E=(E_0,E_1,\ldots ,E_r)$ be an exceptional configuration on $X_r$;
 $r\geq 3$; and let $H\equiv dE_0-m_1E_1-\cdots -m_rE_r$ be an $\bm
 E$-standard class satisfying $\gchi (X_r,\fasm O(H))\geq 3k$ and
 $m_r\geq k-1$, then
 $H$ is non-special and separates $k$-clusters if the following
 necessary conditions are satisfied :  $H.E\geq 2a-1+k$ for all
 $\bm E$-standard isolated curves of genus $a$ ($1\leq a\leq k$)and  $H.E\geq a-1+2k$ for all
 $\bm E$-standard isolated curves of genus $a$ ($\frac{4}{3}k\geq a
 >k$)} 
\begin{prop} Conjecture 2 is equivalent to conjecture 1. \end{prop}
\begin{proof}$\!\!\!\!\!.$ We have already noted in~\ref{minrat}, that for $r\geq
 3$, $H.E_r$ is the 
 minimum value of $H.E$ for all effective rational classes on $X_r$. It
 therefore suffices to
 show that if the stated condition holds on all $\bm E$-standard
 isolated curves of genus  $a$, $1\leq a\leq k$ and $\frac{4}{3}k\geq a
 >k$,  then it holds for all integral curves of such genus. 

An integral curve $F$ on $X_r$ of genus $a\geq 1$ is
standard by~\ref{effst} and $H.F$ is minimum when $F$ is $\bm
E$-semi-standard. If $F$ is $E$-standard, but not isolated, then
letting $j=\mbox{min}\left\{i|1\leq i\leq r \mbox{ and }
  F.E_i=0\right\}$, the class $F\up =F-E_i-\cdots -E_{i+s}$ is an
$\bm E$-standard isolated curve of genus $a>0$
for some $s\geq 0$ and satisfies $F\up.H\leq F.H$. Hence $H.F$ takes its
minimum value on the $\bm E$-standard isolated curves of genus
$a$.\end{proof}
\begin{rem}\be
\item {\em To see that one does have to test in the range $\frac{4}{3}k\geq
  a>k$ when $k\geq 3$, consider $H=13E_0-9E_1-2E_2-\cdots -2E_{18}$ which is
  constructed from the isolated hyperelliptic curve of genus 4,
  $E=6E_0-4E_1-E_2\cdots -E_{18}$. We have $\gchi (\fasm O(H))=3.3$
  and $H.E=8<4-1+2.3$.}  
\item {\em The following example shows that it is necessary to consider
  isolated curves which may not lie on $X_r$. It suffices to
  consider $kC_8$ on $X_8$. For any $k\geq 6$, $\gchi (kC_8)\geq 3k$,
  but $(kC_8).C_8=k<k+1$. However with the definitions introduced in
  (\ref{nd}) this can be detected on the isolated curve $C_9$ on $X_9$.} 
\ee
\end{rem}
\section{Further motivation}
Here we show that the general adjunction theorems imply the conjecture
under much heavier restrictions on $H$.

Recall the
\begin{prop}(see \cite{B-S}, Theorem 2.1) Let $H$ be a divisor
    class on a smooth surface $S$ such that $H-K$ is nef and big and $
  (H-K)^2\geq 4k+1$. Then $H$ separates $k$-clusters unless there
  exists an effective divisor $D$ on $S$ of arithmetic genus $p$ such
  that $H-K-2D$ is effective, 
\begin{equation}\label{B-S} H.D \leq 2p-2+k\quad\mbox{ and }\quad
  2p-2+D^2<H.D <2k+D.K\end{equation}
\end{prop}

Now suppose that H-H holds and let $H$ be a divisor class on $X_r$
such that $\gchi (X_r, \fasm O(H))\geq 3k$, $H.F\geq k-1$ for all
exceptional curves and  $H.F\geq k+1$ for all integral elliptic curves
$F$. Suppose further that $H-K_r$ is nef and big and that  $(H-K_r)^2\geq
4k+1$. Let $D$ be an effective divisor on $X_r$ satisfying (\ref{B-S})
so that $D=E+\sum_{i=1}^s n_iF_i$ is an orthogonal sum of a standard
class $E$ and multiples $n_i>0$ of exceptional curves $F_i$. By
(\ref{B-S}), we find
$$ H.E +\sum_{i=1}^s n_iH.F_i\leq 2(p(E)-n_i(n_i+1)/2-1)-2+k$$
so that $H.E\leq 2p(E)-2$ and $E$ is neither elliptic nor rational,
hence is integral of genus $a\geq 2$. As well, using the
effectiveness of $H-K-2D$, we get $(H-K-2E-2(\sum_{i=1}^s
n_iF_i)).E\geq 0$ so that $(H-K).E>2E ^2$. Now by the last part of (\ref{B-S})
$2k>(H-K).E+(H-K).(\sum_{i=1}^s n_iF_i)\geq 2 E ^2$, so that $k\geq a$
in accordance with the conjecture. 
\section{The not very ample standard classes}
In this section we look at the $\bm E$-standard classes $H$ on $X_r$
with $m_r\geq 1$, $\gchi (X_r,\fasm O(H))\geq 6$ and $H.E< 2a-1+k$ for
some isolated standard class $E$ of genus $a=1$ or $2$ for $k=1,2$;
i.e. the standard classes which have sufficient sections and test
positive on all exceptional curves, but are not base point free or not
very ample because this is not the case for their restriction to an isolated
curve of genus one or two. We will use the list (\ref{isol}) of
isolated, genus two curves $G_i$.

Let $H\equiv dE_0-m_1E_1-\cdots -m_rE_r$ be an $\bm E$-standard class
with $m_r>0$ and $\gchi (X_r,\fasm O(H))\geq 6$.

The only isolated $\bm E$-standard curve of genus one is $C_9=3E_0-E_1-\cdots
-E_9$. To determine the required classes we can suppose that
$r=9$. One easily sees that the standard classes with $H.C_9\leq 1$
and $\gchi(X_r,\fasm o(H))\geq 3$ is  $H\equiv mC_9+E_9$ which has
$\gchi=m$,  and that those with $\gchi (X_r,\fasm O(H))\geq 6$ et $H.C_9\leq 2$, are
$$H\equiv mC_9+E_9;\;\; H\equiv mC_9+2E_9;\;\; H\equiv mC_9+E_8+E_9$$ 
which have $\gchi$, $m+1$, $2m$ and  $2m+1$ respectively.

Let $E$ be an $\bm E$-standard isolated curves of genus $2$
 (\ref{isol}) such that $H.E\leq 4$. If $H.E\leq 3$ then $H-4E$ is 
 efective and negative on $E$ which is an ample divisor! We can thus
 suppose $H.E= 4$. In this case  $H-2E$ is effective.

If $H-3E$ is effective, then  $H=3E+A$
with $1=E.A>A ^2$.  As such, either $H=3E+F$ where $F$ is exceptional
and $E.F=1$, or $H-3E$ is standard and isolated. In the; latter case,
either $H=4E$ or $H=3E+A$ where $A$ is an isolated elliptic curve with
$E.A=1$.(eg. $E=6E_0-2E_1-\cdots -2E_8-E_9-E_10-E_11$, $A=C_9$).  

If $H-3E$ is not effective, then $\overline{\gchi}(H)=5$ and $H=2E+A$
 where $E.A=2$. There are two cases. In the first case, $A$ can be decomposed as a sum
 $A=A_1+A_2$ of reduced and irreducible divisors $A_i$ 
 with $1=E.A_i>A_i ^2\geq -1$ so that the $A_i$
 are exceptional or isolated and elliptic (eg. $E=6E_0-2E_1-\cdots
 -2E_8-E_9-E_{10}-E_{11}$ and 
$$A_1=C_9, A_2=E_{10}\quad \mbox{ or }\quad A_1=E_9,A_2=E_{10})$$
In the second case $H=2E+A$ where $A$ is reduced and
 irreducible, and $4=(E.A)^2>A^2\geq -1$, in which case, $A$ is
 an exceptional curve or an isolated curve of genus $a=$1,2 or 3. One has
 examples with ($E=G_1, A=C_9$), ($E=G_2, A=G_3$) and the pair $E=G_2$
$$A=12E_0-4E_1-\cdots -4E_8-2E_9-2E_{10}-2E_{11}\vspace{5ex}$$

\begin{flushleft}
James E. Alexander\\
Dept. de mathématiques\\
Faculté des Sciences\\
2, bd Lavoisier\\
49045 ANGERS, FRANCE
$\,$\\
jea@univ-angers.fr
\end{flushleft}

\end{document}